\documentclass[12pt,twoside]{article}

\date{}

\usepackage{amssymb}
\usepackage{amsmath}
\usepackage{amscd}
\usepackage{graphicx}
\usepackage{epsfig}
\usepackage{color}
 \newcommand{\R}{\mathbf{\mathbb{R}}}
 \newcommand{\C}{\mathbf{\mathbb{C}}}
 \newcommand{\N}{\mathbf{\mathbb{N}}}

\begin{document}

%\centerline{\bf Journal's Title, Vol. x, 20xx, no. xx, xxx - xxx}

%\centerline{\bf HIKARI Ltd, \ www.m-hikari.com}

%\centerline{\bf http://dx.doi.org/10.12988/}

%\centerline{}

%\centerline{}

\centerline{\Large{\bf On the matrix $pth$ root functions and generalized}}

\centerline{}

\centerline{\Large{\bf Fibonacci sequences}}

\centerline{}

\centerline{\bf {Rajae Ben Taher}}

\centerline{}

\centerline{Group of DEFA - Department of Mathematics and Informatics
    Faculty of Sciences}

\centerline{University of My Ismail}

\centerline{B.P. 4010, Beni M'hamed, Meknes - Morocco}

\centerline{bentaher89@hotmail.fr}

\centerline{}

\centerline{\bf {Youness El Khatabi}}

\centerline{}

\centerline{Group of DEFA - Department of Mathematics and Informatics
    Faculty of Sciences}

\centerline{University of My Ismail}

\centerline{B.P. 4010, Beni M'hamed, Meknes - Morocco}

\centerline{elkhatabi.youness@gmail.com}

\centerline{}

\centerline{\bf {Mustapha Rachidi}}

\centerline{}

\centerline{Group of DEFA - Department of Mathematics and Informatics
    Faculty of Sciences}

\centerline{University of My Ismail}

\centerline{B.P. 4010, Beni M'hamed, Meknes - Morocco}

\centerline{mu.rachidi@hotmail.fr}

\newtheorem{Theorem}{\quad Theorem}[section]

\newtheorem{Definition}[Theorem]{\quad Definition}

\newtheorem{Corollary}[Theorem]{\quad Corollary}

\newtheorem{Lemma}[Theorem]{\quad Lemma}

\newtheorem{Example}[Theorem]{\quad Example}

\newtheorem{Proposition}[Theorem]{\quad Proposition}

\newtheorem{Proof}[Theorem]{\quad Proof}

\newtheorem{remark}[Theorem]{\quad Remark}

\centerline{}

%{\footnotesize Copyright $\copyright$ 20xx Author 1 and Author 2. This article is distributed under %the Creative Commons Attribution License, which permits unrestricted use, distribution, and %reproduction in any medium, provided the original work is properly cited.}

\begin{abstract}
This study is devoted to the polynomial representation of the matrix $p$th root functions. The Fibonacci-H\"{o}rner decomposition of the matrix powers and some techniques arisen from properties of generalized Fibonacci sequences, notably the Binet formula,  serves as a triggering factor to provide explicit formulas for the matrix $p$th roots. Special cases and illustrative numerical examples are given.
\end{abstract}



\section{Introduction}
The $pth$ root of a square matrix occurs in various fields of mathematics, applied sciences, and engineering. For example, this matrix is involved in control and systems theory, matrix differential equations, nonlinear matrix equations, finance and health care. Many methods and techniques have been expanded to provide exact and approximate representations of the matrix $pth$ root (see \cite{abd-bent-khat-rach-2014}, \cite{ben-et-all-2014}, \cite{cross-lancaster-lama-71}, \cite{smith-2003}, and references therein). In this study, we consider the Fibonacci-H\"{o}rner decomposition of the matrix powers (see \cite{abd-bent-rach-2012}, \cite{bent-rach-2008}, \cite{bent-rach-2006} and \cite{bent-rach-2003}) and some techniques based on some properties of generalized Fibonacci sequences (see \cite{dubeau-mrs-1997} and \cite{stanley-1997}),   to provide some explicit formulas of the matrix $p$th roots.
 %%%%
 \par
 %%%%
 Let $A$ be a matrix in $M_d(\C)$, the algebra of $d \times d$ matrices with complex entries $(d\geq 2)$, and $p\geq 2$ a positive integer. Usually {\it a matrix $p$th root of $A$}, is defined as a  matrix $X\in M_d(\C)$ satisfying the equation,
 \begin{equation}\label{eq:X^p=A}
  X^p=A.
 \end{equation}
A matrix $pth$ root can be defined using several definitions of a matrix function of the current literature (see \cite{gantmacher-1960}, \cite{galub-VanLoan-1996}, \cite{horn-johnson-94}, \cite{rinehart-1955} and \cite{verde-star-05}). In general, a matrix $pth$ root may not exist or there may be an infinite number of solutions for (\ref{eq:X^p=A}). In this study, we are particularly interested in the polynomial solutions of Equation (\ref{eq:X^p=A}), when $A$ is nonsingular, in other words the matrix $pth$ roots that are expressible as polynomials in $A$. Such solutions are polynomial functions of a matrix, known as {\it primary matrix functions} (see  \cite{gantmacher-1960}, \cite{horn-johnson-94} and \cite{rinehart-1955}). The function considered here is nothing else but only the complex $pth$ root function $f(z)\equiv z^{1/p}$, which is a multi-valued function. Indeed, for every non-zero complex number
 $z=|z| \exp[i \arg(z)]$ $(-\pi<\arg(z) \leq \pi)$, it is well known that $z$ admits $p$ $pth$ roots given through the use of the functions
 \begin{equation}\label{eq:pth-root-fct}
   f_j(z) = |z|^{1/p} \exp(i[\arg(z)+2\pi j]/p)
   = z^{1/p}\exp(2i\pi j/p), j\in R(p),
 \end{equation}
where  $R(p)=\{0,1,\ldots,p-1\}$. Since $f(z) \equiv z^{1/p}$ is defined on the spectrum of any nonsingular matrix $A$ (see \cite[Ch. 5]{gantmacher-1960}) and giving a choice of $p$ branches for each eigenvalue $\lambda_j$ $(1\leq j\leq l)$ of $A$, many polynomial solutions of Equation (\ref{eq:X^p=A}) may be furnished. To emphasize, the matrix $A$ has precisely $p^s$ matrix $pth$ roots that are primary matrix functions, classified by specifying which branch of the $pth$ root function is taken in the neighborhood of each eigenvalue $\lambda_j$ (see \cite{gantmacher-1960}, \cite{galub-VanLoan-1996}, \cite{horn-johnson-94}, \cite{judith-et-all-2014} and \cite{smith-2003}). In particular, the unique matrix $X$, solution of (\ref{eq:X^p=A}), such that its eigenvalues are in $\{z \in \C \backslash \{0\}: |\arg(z)|<\pi/p\}$ is called {\it the principal matrix $pth$ root of $A$} and will be denoted by $A^{1/p}$. For the sake of simplicity, this primary matrix function will also be denoted by $f(A)=A^{1/p}$.
  \par
Consider a matrix $A \in M_d(\C)$ and a nonzero real parameter $t$ satisfying $\rho(tA)<1$, where $I_d$ is the identity matrix and $\rho(A)$ denote the spectral radius of $A$. The matrix function $g(tA)=(I_d-tA)^{1/p}$ may be defined from the Taylor series expansion $(1-z)^{1/p} = \sum_{n=0}^{\infty} b_{n}z^{n}$, which converges on the open disk $\mathcal{D}(0,1)=\{z\in \C : |z|<1\}$, where $b_{0}=1$ and $b_{n}=(-1)^{n}\frac{\frac{1}{p}(\frac{1}{p}-1)\cdots  (\frac{1}{p}-n+1)}{n!} <0$ for $n \geq 1$. That is to say
  \begin{equation}\label{eq:g(tA)}
    g(tA)=(I_d-tA)^{1/p}=\sum_{n=0}^{\infty} b_{n} t^n A^{n}.
  \end{equation}
  The matrix power series expansion of the adequate function has been used to study the principal matrix $p$th root (see \cite{abd-bent-khat-rach-2014}) and the principal matrix logarithm (see \cite{abd-bent-rach-2012}). One of our main goals is to determine an explicit formula for the principal matrix $pth$ root function $g(tA)$, based mainly on the formula (\ref{eq:g(tA)}) and the Fibonacci-H\"{o}rner decomposition. We highlight that our approach for computing the principal matrix $pth$ root function does not necessarily require the knowledge of the minimal polynomial. Indeed, by employing the characteristic polynomial or any nonzero annihilator polynomial $P(z)$ (of degree $r$), the $nth$ power of $A$ $(n\geq r)$ may be expressed as a linear combination in the Fibonacci-H\"{o}rner system associated to $A$; where the scalar coefficients are the terms of a particular $r$-generalized Fibonacci sequence. By substituting this expression in (\ref{eq:g(tA)}), the Fibonacci-H\"{o}rner decomposition of the principal matrix $p$th root of $I_d-tA$ is obtained. Then the application of the Binet's formula of the $r$-generalized Fibonacci sequence, mentioned above, leads to derive an explicit compact representation of the matrix function $g(tA)$, defined in (\ref{eq:g(tA)}). As a result, by setting $t=1$ and $A=I_d-B$, the principal matrix $pth$ root of $B$ is determined.
  \par
 The remainder of this study is organized as follows. In Section \ref{sect:horner-decomp}, an explicit expression of $(I_d-tA)^{1/p}$ is provided using the Fibonacci-H\"{o}rner decomposition approach. Section \ref{sect: special-cases} is devoted to the presentation of some basic special cases illustrating the method of the preceding section for the computation of the principal matrix $pth$ root. In Section \ref{sect:general-settings} we discuss the polynomial decompositions of {the} primary matrix $pth$ root functions, that satisfy (\ref{eq:X^p=A}), for nonsingular matrices reduced to their Jordan canonical forms. Examples and applications are provided.
 %%%%
\section{Fibonacci-H\"{o}rner decomposition of the principal matrix $p$th root }
\label{sect:horner-decomp}
\subsection{Fibonacci-H\"{o}rner decomposition of $g(tA)=(I_d-tA)^{1/p}$}\label{subsect:fibonacci-horner}
  Let $A$ be in $M_d(\mathbf{\mathbb{C}})$ and a polynomial
$P(z)=z^{r}-a_{0}z^{r-1}-\cdots-a_{r-1}$ $(a_{r-1}\not=0)$ such that
$P(A)=\Theta_{d}$ (zero matrix). The  H\"{o}rner polynomials associated to $P(z)$ are given by $P_0(z)=1$, $P_{j+1}(z)=zP_j(z)-a_{j}$ $(j=0;1;...; r-1)$ and the H\"{o}rner system associated to $A$ is given by $A_0=P_0(A)=I_d$, $A_1=P_1(A)=A-a_0I_d$, $\cdots$,
$A_{r-1}=P_{r-1}(A)=A^{r-1}-a_0A^{r-2}-\cdots-a_{r-2}I_d$. The Fibonacci-H\"{o}rner decomposition of the powers $A^n$ $(n\geq r)$  is given by,
\begin{equation}\label{eq:power_A-u_n}
A^{n}=u_{n}A_{0}+u_{n-1}A_{1}+...+u_{n-r+1}A_{r-1},\; \mbox{ for }\;
n\geq r
\end{equation}
where $u_0=1$, $u_{-1}=\cdots =u_{-r+1}=0$ (see \cite{abd-bent-rach-2012}, \cite{bent-rach-2006} and \cite{bent-rach-2003}).  For
every $n\geq 1$, we show that the term $u_n$ satisfies the linear recursive relation of order $r$ of Fibonacci type,
\begin{equation}\label{eq:r-gfs}
u_{n+1}=a_{0}u_{n}+\cdots+a_{r-1}u_{n-r+1},
\end{equation}
where  $a_0, a_1,\cdots , a_{r-1}$ are specified as the coefficients
of  $\{u_n\}_{n\geq -r+1}$ (see \cite{dubeau-mrs-1997}). With the aid of
Expression (\ref{eq:power_A-u_n}),  we are led to the Fibonacci-H\"{o}rner decomposition of the matrix function $g(tA)$.
 %%%%%%%%%%%%%%%%%%%%%%%%%%%%%%%%%%%%%%THEOREM 1 %%%%%%%%%%%%%%%%%%%
\begin{Theorem}\label{thm : matr-pth-root_Horner}
{\sc  Fibonacci-H\"{o}rner decomposition.}
 Let $A$ be in $M_d(\mathbf{\mathbb{C}})$  and let
$P(z)=z^{r}-a_{0}z^{r-1}-\cdots-a_{r-1}$ $(r\geq 2,
a_{r-1}\not=0)$ be an annihilator polynomial of $A$, i.e.
satisfying $P(A)=\Theta_{d}$. Let $\left\{ A_{s}\right\} _{0\leq
s\leq r-1}$ be the Fibonacci-H\"{o}rner system associated to $A$.
Then, for every $t\in \R \backslash \{0\}$ such that
 $|t|\rho(A)<1$, we
have
\begin{equation}\label{eq:horner-decomp}
g(tA)=(I_{d}-tA)^{1/p}=\sum_{s=0}^{r-1}\varphi_s(t)A_{s}, \mbox{
where  }\varphi_{s}(t) = \sum_{n=s}^{\infty}u_{n-s}b_nt^{n},
\end{equation}
where the $u_n$ are computed from Expression (\ref{eq:r-gfs}).
\end{Theorem}
 %%%%%%%%%%%%%%%%%%%%%%%%%%%%%%%%%%%%%%%%%%%%%%%% PROOF
 \begin{Proof}
 Since
 $|t|\rho(A)<1$, expression  $
(I_d-tA)^{\frac{1}{p}} = \sum_{n=0}^{\infty} b_{n}t^{n}A^n$ shows
that $(1-tA)^{\frac{1}{p}} = Q(t)+H(t)$, where
 $Q(t) =
\sum_{n=0}^{r-1}b_nA^nt^{n}$ and $H(t) = \sum_{n\geq
r}b_nA^{n}t^{n}$. A straightforward computation, using  Expression
(\ref{eq:power_A-u_n}), permits us to derive
$H(t)=\sum_{s=0}^{r-1}(\sum_{n=r}^{\infty}u_{n-s}b_nt^{n})A_{s}$,
where $\left\{ A_{s}\right\} _{0\leq s\leq r-1}$ is the
Fibonacci-H\"{o}rner system associated to $A$. As
Expression (\ref{eq:power_A-u_n}) is still valid for $n\geq 0$, we
show that
$Q(t)=\sum_{n=0}^{r-1}t^{n}b_nA^n=\sum_{n=0}^{r-1}t^{n}b_n\sum_{s=0}^{n}u_{n-s}A_s$,
because $u_{n-s}=0$ when $n<s$. Hence,
$Q(t)=\sum_{s=0}^{r-1}(\sum_{n=s}^{r-1}u_{n-s}b_nt^{n})A_s$.
Finally, we get $(I_{d}-tA)^{1/p}
=\sum_{s=0}^{r-1}(\sum_{n=s}^{\infty}u_{n-s}b_nt^{n})A_s
=\sum_{s=0}^{r-1}\varphi_{s}(t)A_{s}$. The permutation of the finite
sums $ \sum_{s=0}^{r-1}$ and $\sum_{n=r}^{+\infty}$,  follows by the
uniform convergence of the power series in its convergence disc.
Here the power series  $(1-z)^{1/p} = \sum_{n=0}^{\infty}
b_{n}z^{n}$ converges in $D(0;1)$ and, for any  fixed
polynomial $S(z)$, the power series $ \sum_{n=0}^{+\infty} S(n) b_n
z^{n-j} $ converges  also in $D(0;1)$.
 \end{Proof}
 %%%%%%%%%%%%%%%%%%%%%%%%%%%%%%%%%%%%%%%%%%%%%%%%%% Example 1
\begin{Example}\label{examp_matr_order_3_25}{\sc Algebraic matrix of order 3}.
 Let $A$ be in $M_3(\mathbf{\mathbb{C}})$ fulfilling the conditions
of Theorem \ref{thm : matr-pth-root_Horner}, with $P(A)=\Theta_3$,
where $P(z)=z^3-a_0z^2-a_1z-a_2$ with $a_0, a_1, a_2\in
\mathbb{C}$ $(a_2\not= 0)$. The Fibonacci-H\"{o}rner decomposition of
$g(tA)=(I_3-tA)^{1/p}$ is given as $g(tA) =
\sum_{s=0}^{2}\varphi_{s}(t)A_{s}$, with
\begin{equation*}
\varphi_{0}(t)=\sum_{n=0}^{\infty}u_{n}b_nt^{n}\; , \;
\varphi_{1}(t)=\sum_{n=1}^{\infty}u_{n-1}b_n t^{n}\; \mbox{ and } \;
\varphi_{2}(t)=\sum_{n=2}^{\infty}u_{n-2}b_nt^{n}.
\end{equation*}
The sequence $\{u_n\}_{n\geq -2}$ is such that
$u_{-2}=u_{-1}=0$, $u_{0}=1$ and
$u_{n+1}=a_0u_{n}+a_1u_{n-1}+a_2u_{n-2}$.
\end{Example}
  %%%%%%%%%%%%%%%%%%%%%%%%%%%
For every $n\geq r$ a direct computation leads to
\begin{equation*}
  A^{n}=\sum_{s=0}^{r-1}\left(\sum_{j=0}^{s}a_{r-s+j-1}u_{n-j}\right)A^{s},
\end{equation*}
the so-called {\it polynomial decomposition of
$A^n$}. Therefore, the polynomial decomposition of the matrix
function $g(tA)=(I_{d}-tA)^{1/p}$ can also be provided. Indeed, let
$P(z)=z^{r}-a_{0}z^{r-1}-\cdots-a_{r-1}$ $(r \geq 2,
a_{r-1}\not=0)$ be an annihilator polynomial of $A$. Then, for
every $t\in \R \backslash \{0\}$ such that $|t|\rho(A)<1$ we have $g(tA)=(I_{d}-tA)^{1/p}=
\sum_{s=0}^{r-1}\Omega_{s}(t)A^{s}$, where $\displaystyle
\Omega_{s}(t)= b_st^{s}+\sum_{n=r}^{\infty}v_{n}^{(s)} b_nt^{n}$,
for $s=0,1,\cdots,r-1$ with $v_{n}^{(s)}=\sum_{j=0}^{s}a_{r-s+j-1}u_{n-j}$.
 %%%%%%%%%%%%%%%%%%%%%%%%%%%%%%%%%%%%%%%%%%
\begin{remark}\label{thm : remk-similarity}
Let  $A$, $B$ two similar matrices in
$M_d(\mathbf{\mathbb{C}})$, and let $\left\{ A_{s}\right\} _{0\leq
s\leq r-1}$, $\left\{ B_{s}\right\}_{0\leq s\leq r-1}$ be the
Fibonacci-H\"{o}rner systems associated to $A$ and $B$ respectively.
Then, a direct verification shows that $A_s$ and $B_s$ are also
similar. Therefore, if $A$ or $B$
satisfies the conditions of Theorem \ref{thm :
matr-pth-root_Horner}, then these two matrices admit similar
Fibonacci-H\"{o}rner decomposition. Some interesting practical situations could be
studied, throughout similarity of matrices.
\end{remark}
 %%%%%%%%%%%%%%%%%%%%%%%%%%%%%% SUBSECTION
\subsection{Main result}\label{subsect:Binet-method}
Let $A$ be in $M_d(\mathbf{\mathbb{C}})$ such that
$P(A)=\Theta_{d}$, where $P(z)=z^{r}-a_{0}z^{r-1}-\cdots-a_{r-1}$
($r \geq 2$, $a_{r-1}\not=0$). For the sequence $\{u_n\}_{n\geq 0}$ satisfying the
recursive relation of Fibonacci type  (\ref{eq:r-gfs}),  the Binet
formula of its  general term is $\displaystyle  u_n=\sum_{i=1}^{l}\sum_{j=0}^{m_i-1}C_{i,j}n^j\lambda_i^n,$
for all $n\in \mathbb{N}$, where the $\lambda_i$ ($1\leq i\leq l
\leq d$) are the pairwise distinct roots of the polynomial $P(z)$, with algebraic
multiplicities $m_i$ ($\sum_{i=1}^l m_i=r$) (see \cite{dubeau-mrs-1997},  \cite{kelly-peterson-1991} and
\cite{stanley-1997} for example). The coefficients $C_{i,j}$ are obtained by solving the following
linear system of equations,
\begin{equation*}%\label{eq:syst-unk-c(i,j)}
\sum_{i=1}^{l}\sum_{j=0}^{m_{i}-1}C_{i,j}n^j\lambda_i^n=u_n, n=0,
1,...,r-1,
\end{equation*}
 %%%%%%%%%%%%
 An explicit expression of the $C_{i,j}$ can be obtained from a recent work of R. Ben Taher and M. Rachidi \cite{bent-rach-2015-SPMA}.  Indeed, it is shown that,
 %%%
 \begin{equation*}
   \sum_{i=1}^{l}\sum_{j=0}^{m_i-1}C_{i,j}n^j \lambda_i^n = \sum_{i=1}^{l}\sum_{j=0}^{m_i-1}\binom{n+r-1}{j} \gamma_j^{[i]}(\lambda_i)^{n+r-1-j},
 \end{equation*}
 %%%
 with $\gamma_j^{[i]}=\gamma_j^{[i]}(\lambda_1,\ldots,\lambda_l)$ defined by,
 %%%
 \begin{equation*}
 \gamma_j^{[i]}=
 \begin{cases}
   (-1)^{r-m_i}\underset{\sum n_t=m_i-j-1}{\sum}
 \left[ \underset{1\leq t \neq i\leq l}{\prod} \frac{\binom{n_t+m_t-1}{n_t}}{(\lambda_t-\lambda_i)^{n_t+m_t}}\right], & \mbox{if } 1<l\leq d \\
 1 \mbox{ when } j=r-1 \mbox{ and } 0 \mbox{ otherwise},
 & \mbox{if } l=1.
 \end{cases}.
 \end{equation*}
 %%%
 Hence, for a fixed $i=1,\ldots,l$, by aid of the Vandermonde's identity
 \begin{eqnarray*}
 \sum_{j=0}^{m_i-1}C_{i,j}n^j  &=&
 \sum_{j=0}^{m_i-1}\binom{n+d-1}{j} \gamma_j^{[i]}\lambda_i^{r-1-j}
  \\
 % &=& \sum_{j=0}^{m_i-1} \gamma_j^{[i]}\lambda_i^{r-1-j}\sum_{k=0}^{j} \binom{n}{k} \binom{r-1}{j-k} \\
  &=& \sum_{j=0}^{m_i-1} \gamma_j^{[i]}\lambda_i^{r-1-j}
  \sum_{k=0}^{j} \binom{r-1}{j-k}
  \sum_{\ell=0}^{k} \frac{S_{k,\ell}}{k!} n^{\ell},
 \end{eqnarray*}
 where $S_{k,\ell}$ are the \textit{Stirling numbers of the first kind}. Finally, for $i=1,\ldots,l$ and $j=0,\ldots,m_i-1$
 \begin{equation}\label{eq:C-i,j}
  C_{i,j}=\sum_{h=j}^{m_i-1} \lambda_i^{r-1-h} \gamma_h^{[i]} \sum_{k=j}^{h} \binom{r-1}{h-k} \frac{S_{k,j}}{k!}.
  \end{equation}
 %%%%
 \par
 %%%%%%%%%%
 Set $\Delta_1=\{i\, (1
 \leq i \leq l); m_i=1\}$ and $\Delta_2=\{i\, (1\leq i \leq l); 1<
 m_i \leq r\}$.  The Binet formula takes the form,
 %%%%
\begin{equation}\label{eqs : binet_formula_bis}
u_n=\sum_{i\in \Delta_1\cup \Delta_2}\tilde{C}_{i}\lambda_i^n+\sum_{i\in
\Delta_2}\sum_{j=1}^{m_i-1}C_{i,j}n^{j}\lambda_i^{n},
\end{equation}
 %%%
 where $\tilde{C}_{i}=\sum_{i\in \Delta_1}C_{i,0}+
 \sum_{i\in \Delta_2}C_{i,0}$, which may be described as below
 %%%%
 \begin{equation}\label{eq:C_i-tilde}
   \tilde{C}_{i}= \begin{cases}
 \frac{(-1)^{r-1}\lambda_i^{r-1}}
 {\underset{1\leq t \neq i\leq l}{\prod}(\lambda_t-\lambda_i)^{m_t}},
 & \mbox{if } i \in  \Delta_1\\
 \sum_{h=0}^{m_i-1}\lambda_i^{r-1-h}\gamma_h^{[i]}\binom{r-1}{h}, & \mbox{if } i \in  \Delta_2.
 \end{cases}
 \end{equation}
 %%%
 Notice that $\tilde{C}_{i}$ is directly derived from the general expression of the $C_{i,j}$ (\ref{eq:C-i,j}), it will be introduced in the explicit formula of $(I_d-tA)^{1/p}$, provided by the following Theorem.
 In the particular case when $\Delta_2=\emptyset$, the $\lambda_i$  $(1\leq i\leq d)$ are simple, and thus
 $u_n=\sum_{i=1}^{d}\frac{(-1)^{r-1}\lambda_i^{r-1}}
 {\underset{1\leq t \neq i\leq l}{\prod}(\lambda_t-\lambda_i)^{m_t}}\lambda_i^n$. For reason of generality we suppose in the sequel that $\Delta_1 \neq \emptyset$ and $\Delta_2\neq \emptyset$.
 \par
 The main result of this section is presented as follows.
 %%%%%%%%%%%%%%%%%%%%%%%%%%%%%%%%%%%%%%%%%%%%%%%%%%%%%%%%%%% THM
\begin{Theorem}\label{thm_p-th-root_matr_binet}
 Let $A$ be in $M_d(\mathbf{\mathbb{C}})$ such that $P(A)=\Theta_{d}$,
 where $P(z)=z^{r}-a_{0}z^{r-1}-\cdots-a_{r-1}$
$(a_{r-1}\not=0)$. Then, for every $t\in \R \backslash \{0\}$ with
 $|t|\rho(A)<1$, there results $(I_{d}-tA)^{1/p}=\varphi_0(t)I_d+
\sum_{s=1}^{r-1}\left[\Phi_{s}(t)+\Psi_{s}(t)\right]A_{s}$, where $\varphi_0(t)$,
$\Phi_{s}(t)$ and $ \Psi_{s}(t) $ are given by
 %%%
 \begin{equation*}
 \varphi_0(t)=\sum_{i\in \Delta_1 \cup \Delta_2}\tilde{C}_i(1-\lambda_it)^{1/p}+\sum_{i\in \Delta_2}\sum_{j=1}^{m_i-1}C_{i,j} D^j
(1-\lambda_it)^{1/p},
 \end{equation*}
 %%%
 \begin{equation*}%\label{eq:p-th-root-binet-Phi-s}
 \Phi_{s}(t)=
 \sum_{i\in\Delta_1\cup \Delta_2}\frac{\tilde{C}_i}{\lambda_i^{s}}(1-\lambda_{i}t)^{1/p},
 \end{equation*}
 %}
 %%%
 \begin{equation}\label{eq:p-th-root-binet-Psi-s}
 \Psi_{s}(t) = \sum_{i\in \Delta_2}\sum_{j=1}^{m_i-1}
 \sum_{k=0}^j \frac{C_{i,j}}{\lambda_i^{s}} \binom{j}{k}(-s)^{j-k}D^k
 \left((1-\lambda_it)^{1/p}\right).
\end{equation}
Here $D$ is the differential operator $D=t\dfrac{d}{dt}$ (derivation degree operator), $C_{i,j}$ and $\tilde{C}_i$ are respectively given by (\ref{eq:C-i,j}) and (\ref{eq:C_i-tilde}).
\end{Theorem}
  %%%%%%%%%%%%
 \begin{Proof}
 Substitution of the Binet formula (\ref{eqs : binet_formula_bis}) of $u_{n-s}$ in Expression \eqref{eq:horner-decomp} of the $\varphi_{s}(t)$ (taking into
account that $u_{n-s}=0$ if $n<s$), allows us to show that these functions can be expanded under the form,
\begin{equation*}
  \varphi_s(t)=
\begin{cases}
 \sum_{i\in \Delta_1\cup \Delta_2}\Tilde{C}_i\sum_{n=0}^{\infty}b_nt^{n}\lambda_i^{n}+\sum_{i\in
\Delta_2}\sum_{j=1}^{m_i-1}C_{ij}\sum_{n=0}^{\infty} n^{j}b_nt^{n}\lambda_i^{n},\\ \text{ if } s=0, &\\
\sum_{i\in
\Delta_1\cup \Delta_2}\dfrac{\Tilde{C}_i}{\lambda_i^{s}}\sum_{n=s}^{\infty}b_nt^{n}\lambda_i^{n}+\sum_{i\in
\Delta_2}\sum_{j=1}^{m_i-1}\dfrac{C_{ij}}{\lambda_i^{s}}
\Gamma_{ijs}(t),\\ \text{ if } 1 \leq s \leq r-1,&
\end{cases}
\end{equation*}
%}
where
$\Gamma_{ijs}(t)=\sum_{n=s}^{\infty}b_n(n-s)^{j}t^{n}\lambda_i^{n}$.
Since $\displaystyle |t|\rho(A)<1$ we have
 %%%
 \begin{equation*}
   \varphi_0(t)=\sum_{n=0}^{\infty}u_nb_nt^{n}\lambda_i^{n}=\sum_{n=0}^{\infty}
\left(\sum_{i\in \Delta_1\cup \Delta_2}\Tilde{C}_i\lambda_i^n+\sum_{i\in
\Delta_2}\sum_{j=1}^{m_i-1}C_{i,j}n^{j}\lambda_i^{n}\right)b_nt^{n}.
 \end{equation*}
Therefore, we have
\begin{equation*}
\varphi_0(t)=\sum_{i\in \Delta_1 \cup \Delta_2}\Tilde{C}_i(1-
\lambda_it)^{1/p}+\sum_{i\in \Delta_2}\sum_{j=1}^{m_i-1}C_{i,j} D^j
(1-\lambda_it)^{1/p},
\end{equation*}
where $D$ denotes the operator $D=t \dfrac{d}{dt}$.
Second, for $s \geq 1$, we have $\displaystyle
\varphi_s(t)=\sum_{n=s}^{\infty}u_{n-s}b_nt^{n}\lambda_i^{n}=\sum_{n=s}^{\infty}\left(\sum_{i\in
\Delta_1\cup \Delta_2}\Tilde{C}_i\lambda_i^{n-s}+\sum_{i\in
\Delta_2}\sum_{j=1}^{m_i-1}C_{i,j}(n-s)^{j}\lambda_i^{n-s}
\right)b_nt^{n}$. We show that
 %%%
$\varphi_s(t)=\varphi_{s,1}(t)+\varphi_{s,2}(t)$,
where $\varphi_{s,1}(t)=\sum_{n=s}^{\infty}\sum_{i\in \Delta_1 \cup \Delta_2}\Tilde{C}_i\lambda_i^{n-s}b_nt^{n}$ and
 $\varphi_{s,2}(t)=\sum_{n=s}^{\infty}\sum_{i\in
\Delta_2}\sum_{j=1}^{m_i-1}C_{i,j}(n-s)^{j}\lambda_i^{n-s}
b_nt^{n}$. A  direct computation implies that
 %%%
\begin{equation*}
\varphi_{s,1}(t)=R_{s,1}(t)+\sum_{i\in \Delta_1 \cup \Delta_2}
\dfrac{\Tilde{C}_i}{\lambda_i^{s}} (1-\lambda_it)^{1/p},
\end{equation*}
 %%%
where $R_{s,1}(t)=-\sum_{n=0}^{s-1}\sum_{i\in
\Delta_1\cup \Delta_2}\Tilde{C}_i\lambda_i^{n-s}b_nt^{n}$. A similar computation gives
 %%%
 \begin{equation*}
 \varphi_{s,2}(t)= \sum_{n=s}^{\infty}\sum_{i\in
 \Delta_2}\sum_{j=1}^{m_i-1} \sum_{k=0}^j \dfrac{C_{i,j}}{\lambda_i^{s}}
 (^j_k)(-s)^{j-k}n^k\lambda_i^{n} b_nt^{n}= R_{s,2}(t)+
\Omega_s(t),
 \end{equation*}
where  $R_{s,2}$ is the polynomial
 \begin{equation*}
   R_{s,2}(t)=-
 \sum_{n=0}^{s-1}\sum_{i\in \Delta_2}\sum_{j=1}^{m_i-1} \sum_{k=0}^j
 C_{i,j} (^j_k)(-s)^{j-k}n^k\lambda_i^{n-s} b_n t^{n},
 \end{equation*}
 and
 $\Omega_s(t)=\sum_{i\in \Delta_2} \sum_{j=1}^{m_i-1} \sum_{k=0}^j
 \dfrac{C_{i,j}}{\lambda_i^{s}} (^j_k)(-s)^{j-k}
D^k(1-\lambda_it)^{1/p}$. According to the fact that $u_n=0$ for
$n<0$, it is easy to show that
 \begin{eqnarray*}
 \sum_{\ell=1}^{2}R_{s,\ell}(t) &=& -\sum_{n=0}^{s-1}\left[
 \sum_{i\in\Delta_1\cup \Delta_2}\Tilde{C}_i\lambda_i^{n-s}+
 \sum_{i\in\Delta_2}\sum_{j=1}^{m_i-1}C_{i,j}(n-s)^{j}
 \lambda_i^{n-s}\right]b_nt^{n} \\
    &=& -\sum_{n=0}^{s-1} u_{n-s}b_nt^n =0.
 \end{eqnarray*}
 %%%%%%%%%%%%%%%%%%%%%%%%%%%%%%%%%%%%%%%%%%%%%%%%%%%%%
 Expressions of the functions $\varphi_{s,1}$ and $\varphi_{s,2}$ are
derived from the permutations of  the finite sums $\sum_{i\in
\Delta_1 \cup \Delta_2}$ and $\sum_{i\in
\Delta_2}\sum_{j=1}^{m_i-1} \sum_{k=0}^j$ with the infinite sum
$\sum_{n=0}^{\infty}$. Thus the results of the theorem are achieved.
\end{Proof}
  %%%%%%%%%%%%%%%%%%%%%%%%%%%%%%%%%%%%%%%%%%%%%%%%%%%
The family of functions  $\Psi_{s}(t)$  in (\ref{eq:p-th-root-binet-Psi-s}) can be constructed by an induction
process as follows.
%%%%%%%%%%%%%%%%%%%%%%%%%%%%%%%%%%%%%%%%%%%%%%%%%%%%%%%%%%%%%%%%%%Proposition
\begin{Proposition}\label{thm_p-th-Psi(t)}
Under the data of Theorem \ref{thm_p-th-root_matr_binet}, the
functions $\Psi_{s}(t)$ are given by
\begin{equation*}%\label{eq:p-th-root-binet-Lambda-j-2}
 \Psi_{s}(t) = \Psi_{s,0}(t)+
 \sum_{i\in \Delta_2}\sum_{j=1}^{m_i-1}\sum_{k=1}^j
 \dfrac{C_{i,j}}{\lambda_i^{s}}
 \binom{j}{k}(-s)^{j-k}P_{k,i}(t;\lambda_i)\left(1-\lambda_it\right)^{-k+1/p},
\end{equation*}
where $\Psi_{s,0}(t)=\sum_{i\in \Delta_2}\sum_{j=1}^{m_i-1}
 \frac{C_{ij}}{\lambda_i^{s}}(-s)^{j}
 \left(1-\lambda_it\right)^{1/p}$ and $P_{k,i}(t;\lambda_i)$ are the
 polynomials satisfying the equations,
\begin{equation} \label{eq:p-th-root-binet-Lambda-j-P_k,i}
P_{k+1,i}(t;\lambda_i)=t (1-\lambda_it)\dfrac{ d
P_{k,i}(t;\lambda_i)}{dt } +\dfrac{ \lambda_i(pk-1)}{p }t
P_{k,i}(t;\lambda_i),
\end{equation}
with $P_{1,i}(t;\lambda_i)=- \dfrac{ \lambda_it}{p }$.
\end{Proposition}
%%%%%%%%%%%%%%%%%%%%%%%%%%%%%%%%%%%%%%%%%%%%%%%%%%%%%%%%%%%%%%
Formula (\ref{eq:p-th-root-binet-Lambda-j-P_k,i}) is obtained
by a simple induction. We establish also that  functions  $\Psi_{s}(t)$  given by \eqref{eq:p-th-root-binet-Psi-s} in Theorem \ref{thm_p-th-root_matr_binet} may be formulated under the following  compact form,
\begin{equation}\label{eqs : Psi-Operator-formulation}
 \Psi_{s}(t)=\sum_{i\in \Delta_2}W_{s,i}(D) (1- \lambda t)^\frac{1}{p},
 \mbox{   } (s\geq 1),
\end{equation}
 where $W_{0,i}(D)=\sum_{j=0}^{m_i-1}C_{i,j} D^j$ and
 $W_{s,i}(D)= \sum_{j=1}^{m_i-1}\sum_{k=0}^{j} \frac{C_{i,j}}{\lambda^s }
 (^j_k)(-s)^{j-k}D^k$.
 %%%%%%%%%%%%%%%%%%%%%%%
\begin{Example}{\sc Square root of algebraic matrix of order 2.}\label{example-order-2}
 Let $A$ be in $M_r(\mathbf{\mathbb{C}})$, with $P(A)=\Theta_2$,
where $P(z)=M_A(z)=z^2-a_0z-a_1$ is the minimal polynomial of $A$,
with $a_0, a_1 \in \C \backslash \{0\}$. Moreover, we suppose that $A$
satisfies the conditions of Theorem \ref{thm : matr-pth-root_Horner}. We have
\begin{equation*}
  (I_2-tA)^\frac{1}{2}=(\varphi_0(t)-a_0\varphi_1(t))I_2+\varphi_1(t)A,
\end{equation*}
where the  $\varphi_{s}(t) $ $(s=0, 1)$ are given by (\ref{eq:horner-decomp}). Suppose that $A$ admits two distinct eigenvalues $\lambda_1$ and
$\lambda_2$. It's obvious here that $a_0=\lambda_1+\lambda_2$. A
direct computation using techniques and results of Theorem
\ref{thm_p-th-root_matr_binet} implies that
 %%%
 \begin{equation*}
   \varphi_{0}(t) =  \alpha(1- \lambda_1t)^\frac{1}{2} + \beta(1-
\lambda_2t)^\frac{1}{2} \mbox{ and } \varphi_{1}(t) =
  \frac{\alpha }{\lambda_1} (1-\lambda_1t)^\frac{1}{2} +
 \frac{\beta }{\lambda_2}(1-\lambda_2t)^\frac{1}{2},
 \end{equation*}
 %%%
 where $\alpha = \frac{\lambda_1
}{\lambda_1-\lambda_2}$ and $\beta = \frac{\lambda_2
}{\lambda_2-\lambda_1}$.  Particularly, if $A=\left(\begin{array}{cc}
e  & g \\
f & h
\end{array}
\right) $ satisfies the conditions of  Theorem \ref{thm : matr-pth-root_Horner}, we obtain
 \begin{equation*}
   (I_d-tA)^\frac{1}{2}
 =\left(\begin{array}{cc}
\varphi_0(t)+(e-a_0)\varphi_1(t)  & g \varphi_{1}(t)\\
f\varphi_{1}(t) & \varphi_0(t)+(h-a_0)\varphi_1(t)
\end{array}
\right).
 \end{equation*}
 %%%
Therefore, the square root of the matrix $B=I_r-A$ is given by
 \begin{equation*}
   B^\frac{1}{2}
 =\left(\begin{array}{cc}
\varphi_0(1)+(e-a_0)\varphi_1(1)  & g \varphi_{1}(1)\\
f\varphi_{1}(1) & \varphi_0(1)+(h-a_0)\varphi_1(1)
\end{array}
\right).
 \end{equation*}
 %%%
As a numerical illustration, if we consider  the matrix $B=\left(\begin{array}{cc}
 \frac{1 }{6} & -1 \\
 \frac{1 }{6}  & 1
\end{array}
\right) $, we can check that  $A=I_2-B=\left(\begin{array}{cc}
 \frac{5 }{6} & 1 \\
 \frac{-1 }{6}  & 0
\end{array}
\right) $ and the root of $P_A(z)=z^2- \frac{5 }{6} z+ \frac{1}{6} $ are $\lambda_1= \frac{1}{2} $, $\lambda_2= \frac{1 }{3} $. Therefore, a direct computation shows that the square root of the matrix $B$ is  $ B^\frac{1}{2}
 =\sqrt{2} \left(\begin{array}{cc}
 \frac{3}{2} -\frac{2\sqrt{3} }{3} & 3-2\sqrt{3} \\
 \frac{- 1 }{2} +\frac{\sqrt{3} }{3} &  -1+\sqrt{3}
\end{array}
\right)$.
\end{Example}
 %%%%%%%%%%%%%%%%%%%%%%SECTION 4%%%%%%%%%%%%%%%%%%%%		
\section{Special cases}\label{sect: special-cases} We are interested here  in
the principal matrix $p$th root of a matrix  $A\in M_d(\mathbb{C})$,
whose annihilating-polynomial assumes one of the following forms:
$P(z)=(z-\lambda)^r$ , $P(z)=(z-\mu)(z-\lambda)^{r-1}$   and
$P(z)=(z-\lambda)^{m_1}(z-\mu)^{m_2}$ ($m_1+m_2=r$), where $r \geq 2$. In this subsection we suppose that the matrix
 $A\in M_d(\mathbb{C})$ satisfies the conditions of Theorem
 \ref{thm_p-th-root_matr_binet},  and
 $P(A)=\Theta_d$.
  %%%%%%%%%%%%%%%%%%%%%%%%%%%%%%%%% SUBSECTION 4.1
 \subsection*{Case $P(z)=(z-\lambda)^r$}
 The recurrence relation associated to the
sequence $\{u_n\}_{n\geq 0}$ is defined by
%\begin{equation}\label{equ:rec_rel_z-lambda}
\begin{equation*}
  u_{n+1}=\sum_{s=0}^{r-1}a_{s}u_{n-s}\quad ; \quad
a_{s}=-\binom{r}{s+1}(-\lambda)^{s+1},
\end{equation*}
 %%%
and the  customary initial conditions are $u_0=1$, $u_{n}=0$,  for
$n\leq -1$. We notice here that $\Delta_1=\emptyset $ and
$\Delta_2=\{\lambda \}$. Consequently, the Binet formula (\ref{eqs : binet_formula_bis}), provided
$u_n=\left(C_{0}+\sum_{j=1}^{r-1}C_{j}n^{j}\right)\lambda^{n}$, with
$C_{0}=1$.
 %%%%
\par
 According to the formula (\ref{eq:C-i,j}), the expression of the $C_j$ $(j=0,\ldots,r-1)$ takes the form
\begin{equation*}%\label{eq:C-j}
  C_j=\sum_{h=j}^{r-1} \binom{r-1}{h} \frac{S_{h,j}}{h!}.
\end{equation*}
 %%%%
Therefore, by Theorem \ref{thm_p-th-root_matr_binet}, we
obtain the  result.
%%%%%%%%%%%%%%%%%%%%%%%%%%%%%%% PROPOSTION
\begin{Proposition}\label{thm : prop_R(z)=(z-lambda)-r}
Let $A$ be in $M_d(\mathbb{C})$ satisfying the data of
Theorem  \ref{thm_p-th-root_matr_binet} with $P(A)=\Theta_d$, where
$P(z)=(z-\lambda)^r$ $(r\geq 2)$. Then, for every $t\in\R \backslash \{0\}$, with
 $|t|\rho(A)<1$,
$(I_{d}-tA)^{1/p}=\sum_{s=0}^{r-1}\varphi_{s}(t)A_{s}$, where the
$\varphi_{s}(t)$ are given by $\displaystyle \varphi_{0}(t)=(1-
\lambda t)^\frac{1}{p}+\sum_{j=1}^{r-1}C_jD^j(1- \lambda
t)^\frac{1}{p}$
 and $ \varphi_{s}(t)=\frac{1}{\lambda^s}(1- \lambda t)^\frac{1}{p}
 + \sum_{j=1}^{r-1}\sum_{k=0}^{j}
\frac{C_j}{\lambda^s }(-s)^{j-k}(^j_k)D^k(1- \lambda t)^\frac{1}{p}$
for $s\geq 1$ (recall that $D= t\frac{d}{dt}$).
  \end{Proposition}
 %%%%%%%%%%%%%%%%%%%%%%%%%%%%%%%%%%%%%%%%%%%%% COROLLARY
 As in Expression \eqref{eqs : Psi-Operator-formulation}, the preceding results
 may take the following compact form.
 \begin{Corollary}\label{thm : coro-prop_R(z)=(z-lambda)-r}
Let $D= t\frac{d}{dt}$ be the derivation degree operator. Under the data of
Proposition \ref{thm : prop_R(z)=(z-lambda)-r},
$(I_{d}-tA)^{1/p}=\sum_{s=0}^{r-1}\varphi_{s}(t)A_{s}$. The
$\varphi_{s}(t)$ are given by
 %%%
 \begin{equation*}
   \varphi_{0}(t)=W_0(D)(1- \lambda t)^\frac{1}{p} \mbox{ and }
\varphi_{s}(t)=W_s(D)(1- \lambda t)^\frac{1}{p},  \mbox{ for  }
s\geq 1,
 \end{equation*}
 %%%
and the operators $W_0(D)$ and $W_1(D)$ are given by
 \begin{equation*}
   W_0(D)= 1\!\!1+\sum_{j=1}^{r-1}C_jD^j \mbox{ and }
 W_s(D)=
 \frac{1}{\lambda^s}1\!\!1+\sum_{j=1}^{r-1}\sum_{k=0}^{j}
 \frac{C_j}{\lambda^s } (-s)^{j-k}(^j_k)D^k .
 \end{equation*}
 %%%
We recall that $1\!\!1$ is the identity operator.
 \end{Corollary}
 %%%%%%%%%%%%%%%%%%%%%%%%%%%%%%%%%EXAMPLE
    \begin{Example}
    Compute the principal square root of the matrix,
    \begin{equation*}
      B= \left[ \begin{array}{ccc}
 \frac{10}{6} & -\frac{2}{3}  &  -\frac{1}{3} \\
 \\
 \frac{7}{12} &  \frac{1}{6}  &  -\frac{1}{6} \\
 \\
 \frac{35}{12} &   -\frac{5}{3} &  -\frac{1}{3} \\
  \end{array}\right].
    \end{equation*}
    %%%
  Let consider the matrix $A=I_3-B$, the minimal polynomial of $A$
  takes the form $P(z)=M_A(z)=(z-\frac{1}{2})^2$. Therefore, applying the Proposition
  \ref{thm : prop_R(z)=(z-lambda)-r}  we obtain
  \begin{equation*}
    (I_3-A)^{\frac{1}{2}}=\varphi_0I_3+\varphi_1(A-a_0I_3),
  \end{equation*}
  %%%
  where $\varphi_0=\varphi_0(1)=\frac{\sqrt{2}}{2}+
  C_{1}D(1-\frac{1}{2}t)^{\frac{1}{2}}_{|_{t=1}}$ and
  \begin{equation*}
    \varphi_1=\varphi_1(1)=\sqrt{2}+2C_1\sum_{k=0}^1 (-1)^{1-k}
  D^{k}(1-\frac{1}{2}t)^{\frac{1}{2}}_{|_{t=1}}
  = \sqrt{2}-\frac{3}{2}\sqrt{2}C_1.
  \end{equation*}
  %%%
  By a straightforward computation, we show that $a_0=1$ and $C_1=1$.
 Thus, we obtain $\varphi_0=\frac{\sqrt{2}}{4}$ and
 $\varphi_1=-\frac{\sqrt{2}}{2}$.  Consequently, $B^{\frac{1}{2}}=(I_3-A)^{\frac{1}{2}}=\varphi_0I_3+\varphi_1(A-I_3)
=\varphi_0 I_3 - \varphi_1 B$, which allows to compute the entries
of the principal matrix square root of $B$.
     \end{Example}
 %%%
 \begin{remark}
   In the precedent example, any annihilator polynomial of $A=I_3-B$ could be used for computing the square root of the matrix $B$. For instance, by considering the characteristic polynomial of $A$ the expression of the principal matrix square root of $B$ will take the form:
   \begin{equation*}
     B^{1/2}=\tilde{\varphi}_0 I_3+\tilde{\varphi}_1 B+\tilde{\varphi}_2 B^2,
   \end{equation*}
   %%%
  where $\tilde{\varphi}_0=\frac{5\sqrt{2}}{16}$, $\tilde{\varphi}_1=\frac{\sqrt{2}}{4}$ and $\tilde{\varphi}_2=\frac{\sqrt{2}}{2}$. As has been stated, the latter expression permits to obtain the same matrix that was obtained using the minimal polynomial of $A$.
 \end{remark}
  %%%%%%%%%%%%%%%%%%%%%%%%%%%%%%%%%%%%%%%%%%%%%%%%%%%%%%%%%%%%%%
 \subsection*{Case $P(z)=(z-\mu)(z-\lambda)^{r-1}$} The sequence
 $\{u_n\}_{n\geq 0}$ is defined by $u_{n+1}=\sum_{s=0}^{r-1}a_{s}u_{n-s}$
of coefficients
 \begin{equation*}
   a_{r-1}=\mu(-\lambda)^{r-1} \mbox{ and  }
 a_s=\left[(^{r-1}_{s+1})\lambda+(^{r-1}_{s})\mu \right](-\lambda)^s,
 \mbox{ for } s=0,\dots,r-2,
 \end{equation*}
 %%%
and  initial data $u_0=1$, $u_{n}=0$ for $n\leq -1$. Since
$\Delta_1=\{\mu\} $ and
 $\Delta_2=\{\lambda\}$, the Binet formula, Expression
(\ref{eqs : binet_formula_bis})
 yield $u_n= C_{\mu}\mu^n+\left(\sum_{j=0}^{r-2}n^{j}C_{\lambda,j}\right)
 \lambda^{n}$.  Here the coefficients $C_{\mu}$ and $C_{\lambda,j}$ are
 obtained using the formula (\ref{eq:C-i,j}). Consequently, by
 Theorem \ref{thm_p-th-root_matr_binet}, we obtain the  result.
  %%%%%%%%%%%%%%%%%%%%%%%%%%%%%%%%%%%%%%%%%%%%%%PROPOSITION
\begin{Proposition}\label{thm : prop_R(z)=(z-mu)(z-lambda)-(r-1)}
Let $A$ be in $M_d(\mathbb{C})$ satisfying the
conditions of Theorem  \ref{thm_p-th-root_matr_binet} with
$P(A)=\Theta_d$,
 where $P(z)=(z-\mu)(z-\lambda)^{r-1}$ $(r\geq 2)$. Then, for every $t\in \R \backslash \{0\}$
 with  $|t|\rho(A)<1$, there results
$(I_{d}-tA)^{1/p}=\sum_{s=0}^{r-1}\varphi_{s}(t)A_{s}$, where the
$\varphi_{s}(t)$ are given by $
 \varphi_{0}(t)= C_{\mu} (1- \mu t)^\frac{1}{p}+
 \sum_{j=0}^{r-2}C_{\lambda,j}D^j(1- \lambda t)^\frac{1}{p}$, and
 \begin{equation*}
   \varphi_{s}(t)=\frac{C_{\mu}}{\mu^s } (1- \mu t)^\frac{1}{p}+
 \sum_{j=0}^{r-2}\sum_{k=0}^{j}\frac{C_{\lambda,j}}{\lambda^s}
 (-s)^{j-k}(^j_k)D^k(1- \lambda t)^\frac{1}{p},
 \end{equation*}
 %%%
for $s\geq 1$, where $D= t\frac{d}{dt}$ and $D^0=1\!\!1$ the
identity operator.
 \end{Proposition}
  %%%%%%%%%%%%%%%%%%%%%%%%%%%%%%%%%%%%%%%%%%%%%%%%
Similarly to Expression \eqref{eqs : Psi-Operator-formulation} and
Corollary \ref{thm : coro-prop_R(z)=(z-lambda)-r}, the functions
$\varphi_{s}(t)$ can be expressed as follows $ \varphi_{0}(t)=
C_{\mu} (1- \mu t)^\frac{1}{p}+
 W_{0,\lambda}(D)(1- \lambda t)^\frac{1}{p}$ and $\varphi_{s}(t)=\frac{C_{\mu}}{\mu^s }(1- \mu t)^\frac{1}{p}+
  W_{s,\lambda}(D)(1- \lambda t)^\frac{1}{p}$, for $s\geq 1$,
  where $W_{0,\lambda}(D)=\sum_{j=0}^{r-2}C_{\lambda,j}D^j$ and

 \begin{equation*}
   W_{s,\lambda}(D)=\sum_{j=0}^{r-2}\sum_{k=0}^{j}\frac{C_{\lambda,j}}{\lambda^s}
 (-s)^{j-k}(^j_k)D^k=\sum_{j=0}^{r-2}\frac{C_{\lambda,j}}{\lambda^s
 }(D-s)^k, \mbox{ for } s\geq 1.
 \end{equation*}
  %%%%%%%%%%%%%%%%%%%%%%%%%%%%%%%%%%%%%%%%%%%%%%%%%%%%%%%%EXAMPLE
  \begin{Example}
 Compute the principal cubic root of the matrix $B$ such that,
 \begin{equation*}
   B=\left[ \begin{array}{ccc}
 \frac{3}{4} & 1  & -1  \\
 \\
 \frac{1}{24} & \frac{1}{2} & -\frac{1}{6} \\
 \\
  \frac{5}{48} & \frac{1}{4} & \frac{1}{12} \\
  \end{array}\right].
 \end{equation*}
 %%%
 For $A=I_3-B$, we have
 $P_A(z)=M_A(z)=(z-\lambda)^2(z-\mu)$, with $\lambda=\frac{1}{2}$ and
 $\mu=\frac{2}{3}$. Using Proposition \ref{thm : prop_R(z)=(z-mu)(z-lambda)-(r-1)}, we derive
 \begin{eqnarray*}
 B^{\frac{1}{3}}&=&(I_3-A)^{\frac{1}{3}}=\varphi_0 I_3+\varphi_1 A_1+
 \varphi_2 A_2 \\
 &=&\varphi_0 I_3+\varphi_1(A-a_0 I_3)+\varphi_2(A^2-a_0 A-a_1 I_3)\\
 &=& \left[\varphi_0+(1-a_0)\varphi_1+(1-a_0-a_1)\varphi_2\right]I_3+
 \left[(a_0-2)\varphi_2-\varphi_1 \right]B+\varphi_2 B^2,
 \end{eqnarray*}
 where $a_0=\frac{5}{3}$, $a_1=-\frac{11}{12}$ and the formulas of
 $\varphi_i=\varphi_i(1)$, \textup{(i=0,1,2)} are given by
 Proposition \ref{thm : prop_R(z)=(z-mu)(z-lambda)-(r-1)}. Hence, $\varphi_0=C_{\mu}(1-\mu)^{\frac{1}{3}}+\sum_{j=0}^{1}C_{\lambda,j}
 D_{|_{t=1}}^{j}(1-\lambda t)^{\frac{1}{3}}$,
 \begin{equation*}
   \varphi_1=\frac{C_{\mu}}{\mu}(1-\mu)^{\frac{1}{3}}+\sum_{j=0}^{1}
 \sum_{k=0}^{j} \frac{C_{\lambda,j}}{\lambda}(-1)^{j-k}(^j_k)D_{|_{t=1}}^{k}
  (1-\lambda t)^{\frac{1}{3}} \mbox{ and }
 \end{equation*}
 %%%
 \begin{equation*}
   \varphi_2=\frac{C_{\mu}}{\mu^2}(1-\mu)^{\frac{1}{3}}+\sum_{j=0}^{1}
 \sum_{k=0}^{j} \frac{C_{\lambda,j}}{\lambda^2}(-2)^{j-k}(^j_k)D_{|_{t=1}}^{k}
  (1-\lambda t)^{\frac{1}{3}}.
 \end{equation*}
  %%%
 It follows from (\ref{eq:C-i,j}) that
 $C_{\mu}=16$, $C_{\lambda,0}=-15$ and $C_{\lambda,1}=-3$. Thus, by a
 direct computation, we derive
 $\varphi_0=16(\frac{1}{3})^{\frac{1}{3}}-
 15(\frac{1}{2})^{\frac{1}{3}}+2^{-\frac{1}{3}}$,
 $\varphi_1=24(\frac{1}{3})^{\frac{1}{3}}-24(\frac{1}{2})^{\frac{1}{3}}+
 2^{\frac{2}{3}}$ and $\varphi_2=36(\frac{1}{3})^{\frac{1}{3}}-
 36(\frac{1}{2})^{\frac{1}{3}}+ 2^{\frac{5}{3}}$.  Therefore, $B^{\frac{1}{3}}=\alpha I_3+\beta B+\varphi_2 B^{2}$,  where
$\alpha=9(\frac{1}{3})^{\frac{1}{3}} - 8(\frac{1}{2})^{\frac{1}{3}}
 + \frac{2^{\frac{2}{3}}}{3}$ and $\beta=36(\frac{1}{2})^{\frac{1}{3}} -
 36(\frac{1}{3})^{\frac{1}{3}} - \frac{5}{3}2^{\frac{2}{3}}$.
  The entries of the principal matrix $3th$ root can be derived easily.
  \end{Example}
 %%%%%%%%%%%%%%%%%%%%%%%%%%%%%%%%%%%%%%%%%%%%%%%%%%%% SUBSECTION 4.3
\subsection*{Case $P(z)=(z-\lambda_1)^{m_1}(z-\lambda_2)^{m_2}$}
 Let $\{u_n\}_{n\geq 0}$ be the associated recursive
 sequence, with customary initial conditions $u_0=1$, $u_{n}=0$ for
 $n\leq -1$. From the Binet formula, Expressions (\ref{eqs : binet_formula_bis}) yields $u_n=
\left(\sum_{j=0}^{m_1-1}C_{1,j}n^{j}\right)\lambda_1^{n}+
\left(\sum_{j=0}^{m_2-1}C_{2,j}n^{j}\right)\lambda_2^{n}$. Easily , we show that $\Delta_1=\emptyset$ and $\Delta_2=\{\lambda_1,
\lambda_2\}$. Therefore,  applying Theorem \ref{thm_p-th-root_matr_binet}, we obtain the  result.
%%%%%%%%%%%%%%%%%%%%%%%%%%%%%%%%%%%%%%%%%%%%%%%%%%%%Proposition
\begin{Proposition}\label{thm : prop_R(z)=(z-mu)^m1(z-lambda)^(m_2)}
Let $m_1,m_2 \in \mathbb{N}\backslash \{0;1\}$ and $A$ in
$M_d(\mathbb{C})$ satisfying the conditions of Theorem
\ref{thm_p-th-root_matr_binet}  such that $P(A)=\Theta_d$, for
$P(z)=(z-\lambda_1)^{m_1}(z-\lambda_2)^{m_2}$.
Then, for every $t\in \R \backslash \{0\}$ with
 $|t|\rho(A)<1$, we have
$(I_{d}-tA)^{1/p}=\sum_{s=0}^{m_1+m_2-1}\varphi_{s}(t)A_{s}$. The
functions $\varphi_{s}(t)$ are given by $ \displaystyle
\varphi_{0}(t)= \sum_{i=1}^2 \sum_{j=0}^{m_i-1}C_{i,j}
 D^j(1-\lambda_i t)^\frac{1}{p}$
 and
 \begin{equation*}
 \varphi_{s}(t)= \sum_{i=1}^2\frac{C_{i,0}}{\lambda_{i}^s}
 (1-\lambda_{i}t)^{1/p}+\sum_{i=1}^2\sum_{j=1}^{m_i-1}
 \sum_{k=0}^j \frac{C_{i,j}}{\lambda_i^{s}} \binom{j}{k}(-s)^{j-k}D^k
 \left((1-\lambda_it)^{1/p}\right),
 \end{equation*}
 %%%
for $s\geq 1$, where $D= t\frac{d}{dt}$.
\end{Proposition}
  %%%%%%%%%%%%%%%%%%%%%%%%%%%%%%%%%%%%%%%%%%%%
As in  Expression \eqref{eqs : Psi-Operator-formulation} and
Corollary \ref{thm : coro-prop_R(z)=(z-lambda)-r}, the compact
 expressions of the functions  $\varphi_{s}(t)$ can be
given as follows,
 \begin{equation*}
 \varphi_{0}(t)= \sum_{i=1}^2 W_{i,0}(D)(1- \lambda_i t)^\frac{1}{p}
 \mbox{ and }
 \varphi_{s}(t)=\sum_{i=1}^2 W_{i,s}(D)(1- \lambda_i t)^\frac{1}{p}.
 \end{equation*}
 %%%
 where
 $W_{i,0}(D)=\sum_{j=0}^{m_i-1}C_{i,j} D^j$ and
 %%%
 \begin{equation*}
 W_{i,s}(D)=\frac{C_{i,0}}{\lambda_{i}^s}1\!\!1+
 \sum_{j=1}^{m_i-1} \sum_{k=0}^j \frac{C_{i,j}}{\lambda_i^{s}}
 \binom{j}{k}(-s)^{j-k}D^k .
 \end{equation*}
 %%%
 Note that the result of Proposition
   \ref{thm : prop_R(z)=(z-mu)^m1(z-lambda)^(m_2)} is nothing else but the
   superposition of two cases of the formula given in Proposition
   \ref{thm : prop_R(z)=(z-lambda)-r}.
  %%%%%%%%%%%%%%%%%%%%%%%%%%%%%%%%%EXAMPLE
  \begin{Example}
  Let compute the principal $4th$ root of the matrix,
  \begin{equation*}
  B=\left[ \begin{array}{cccc}
 \frac{3}{4} & 7 & -1 & -\frac{3}{2} \\
 \\
 0 & 0 & 0 &  \frac{1}{8} \\
 \\
 \frac{1}{8} & \frac{15}{4} & 0 & -\frac{7}{8} \\
 \\
 0 & -1 & 0 &  \frac{3}{4}
  \end{array}\right].
  \end{equation*}
  %%%
For $B=I_r-B$ we have $P_A(z)=(z-\frac{1}{2})^2(z-\frac{3}{4})^2$,
thus application of the Proposition \ref{thm : prop_R(z)=(z-mu)^m1(z-lambda)^(m_2)}
leads to the formula,
 \begin{eqnarray*}
 (I_4-A)^{\frac{1}{4}} &=&\left[\varphi_0-a_0\varphi_1-a_1\varphi_2-a_2\varphi_3 \right] I_4+
 \left[\varphi_1-a_0\varphi_2-a_1\varphi_3 \right]A \\
 &&+
 \left[\varphi_2-a_0\varphi_3 \right]A^2+\varphi_3 A^3,
 \end{eqnarray*}
 where the $a_i$ $(i=0,1,2,3)$ are calculated from the
 coefficients of the polynomial
 $P(z)=(z-\frac{1}{2})(z-\frac{3}{4})$.
 Therefore, $a_0=\frac{5}{2}$,
 $a_1=-\frac{37}{16}$, $a_2=\frac{15}{16}$ and $a_3=-\frac{9}{64}$.
 The formulas of the coefficients $\varphi_s$, are given by
 $\varphi_0=\sum_{i=1}^2 \sum_{j=0}^{1}C_{i,j}
 D_{|_{t=1}}^j(1-\lambda_i t)^\frac{1}{4}$ and
 $\varphi_s=\sum_{i=1}^2\frac{C_{i,0}}{\lambda_{i}^s}
 (1-\lambda_{i})^{\frac{1}{4}}+\sum_{i=1}^2
 \sum_{k=0}^1 \frac{C_{i,1}}{\lambda_i^{s}} (-s)^{1-k}
 D_{|_{t=1}}^k \left((1-\lambda_it)^{1/p}\right)$ for $s=1,2,3$,
 here $\displaystyle \lambda_1=\frac{1}{2}$ and $\lambda_2=\frac{3}{4}$.
 The following coefficients are computed by applying Formula (\ref{eq:C-i,j})
 $C_{1,0}=28$, $C_{1,1}=4$, $C_{2,0}=-27$ and $C_{2,1}=9$.
 Therefore, we have the numerical result,
 $\varphi_0=28(\frac{1}{2})^{\frac{1}{4}}-27(\frac{1}{4})^{\frac{1}{4}}
 -\frac{1}{2}(\frac{1}{2})^{-\frac{3}{4}}-
 \frac{27}{16}(\frac{1}{4})^{-\frac{3}{4}}$,
 $\varphi_1=48(\frac{1}{2})^{\frac{1}{4}}-48(\frac{1}{4})^{\frac{1}{4}}
 -(\frac{1}{2})^{-\frac{3}{4}}-\frac{9}{4}(\frac{1}{4})^{-\frac{3}{4}}$,
 $\varphi_2=80(\frac{1}{2})^{\frac{1}{4}}-80(\frac{1}{4})^{\frac{1}{4}}
 -2(\frac{1}{2})^{-\frac{3}{4}}-3(\frac{1}{4})^{-\frac{3}{4}}$
 and
 $\varphi_3=128(\frac{1}{2})^{\frac{1}{4}}-128(\frac{1}{4})^{\frac{1}{4}}
 -4(\frac{1}{2})^{-\frac{3}{4}}-4(\frac{1}{4})^{-\frac{3}{4}}$.
 Hence,  the principal matrix $4-th$ root of $B$ is given as follows
 \begin{equation*}
   B^{\frac{1}{4}}= \alpha I_4+\beta (I_4-B)+\gamma
 (I_4-B)^2+\varphi_3 (I_4-B)^3,
 \end{equation*}
 where $\alpha=- 27(\frac{1}{2})^{\frac{1}{4}}+28(\frac{1}{4})^{\frac{1}{4}}
 +\frac{9}{8}(\frac{1}{2})^{-\frac{3}{4}}+\frac{3}{4}(\frac{1}{4})^{-\frac{3}{4}}$,
 $\beta=144(\frac{1}{2})^{\frac{1}{4}}-144(\frac{1}{4})^{\frac{1}{4}}
 -\frac{21}{4}(\frac{1}{2})^{-\frac{3}{4}}-4(\frac{1}{4})^{-\frac{3}{4}}$ and
 $\gamma=-240(\frac{1}{2})^{\frac{1}{4}}+240(\frac{1}{4})^{\frac{1}{4}}
 +8(\frac{1}{2})^{-\frac{3}{4}}+7(\frac{1}{4})^{-\frac{3}{4}}$. Finally, the entries of the matrix $ B^{\frac{1}{4}}$
 can  be obtained by a direct calculation.
 \end{Example}
 %%%%%%%%%%%%%%%%%%%%%%%%%%%%%%%%%%%%%%%%%%%%%%%%%%%%%SECTION 4
 \section{General settings}\label{sect:general-settings}
 %%%%%
In this section we make use of the results performed in sections \ref{sect:horner-decomp} and \ref{sect: special-cases}, combined with \cite[Theorem 2]{ben-rach-2002}, to determine the polynomial decompositions of all primary matrix $pth$ root functions of a non-singular matrix.
 \par
 The following result is a direct consequence of Proposition \ref{thm : prop_R(z)=(z-lambda)-r}.
 %%%%%%%%%%%% Corollary
 \begin{Proposition}\label{thm : prop-f_j}
   Let $A$ be in $M_d(\C)$ with minimal polynomial $M_A(z)=(z-\lambda)^m$ and $j \in R(p)=\{0,1,\ldots,p-1\}$. Then, for every non-zero parameter $t$ such that $|t\lambda|<1$
 \begin{equation*}
 f_j(I_d-t A)=\exp(2i \pi j/p)\sum_{s=0}^{m-1}
 \varphi_s(t) A_s,
 \end{equation*}
 %%%
 where $f_j(z)$  $(j\in R(p))$ are given as in (\ref{eq:pth-root-fct}), $\displaystyle \varphi_{0}(t)=(1-
 \lambda t)^\frac{1}{p}+\sum_{j=1}^{m-1}C_jD^j(1- \lambda
 t)^\frac{1}{p}$,
 $ \varphi_{s}(t)=\frac{1}{\lambda^s}(1- \lambda t)^\frac{1}{p}
 + \sum_{j=1}^{m-1}\sum_{k=0}^{j}
 \frac{C_j}{\lambda^s }(-s)^{j-k}(^j_k)D^k(1- \lambda t)^\frac{1}{p}$
 $(\textup{for } s\geq 1)$ and $C_j=\sum_{h=j}^{m-1} \binom{m-1}{h} \frac{S_{h,j}}{h!}$.
 \end{Proposition}
 %%%%%%%%%%%%%%%%%%%%%%%%%%%%%%%%%%%%%%%%%%%%%%%%%%%
 Consider $A \equiv \bigoplus_{k=1}^{l} M_k$ (where $M_k = M_k(\lambda_k)$)
 with minimal polynomial $M_A(z)=\prod_{i=1}^{l}(z-\lambda_i)^{m_i}$ $(\sum_{i=1}^{l} m_i \leq d)$,
 and a non-zero real parameter $t$ such that $|t|\sigma(A)<1$. It is known that the $M_k$ are of minimal polynomial $(z-\lambda_k)^{m_k}$ and characteristic polynomial $(z-\lambda_k)^{d_k}$, where $d_k$ is the sum of sizes of all Jordan blocks associated to $\lambda_k$.
 It follows from the previous proposition that,
 %%%%%%%%%%%%%%%%%
 \begin{eqnarray*}
 f_{[j]}(I_d-tA) &=& \bigoplus_{k=1}^{l} f_{j_k}(I_{d_k}-t M_k) \\
    &=& \bigoplus_{k=1}^l \exp(2i \pi j_k/p)
 \left[\varphi_{(k,0)}(t)I_{d_k}+\sum_{\tau=1}^{m_k-1} \varphi_{(k,\tau)}(t) M_{(k,\tau)}\right],
 \end{eqnarray*}
  %%%%%%%%%%%%%%%
 where $[j]\equiv (j_1,\ldots,j_l)$ $(j_k \in R(p))$ and $M_{(k,\tau)}=M_k^{\tau}-\sum_{\ell=0}^{\tau-1} a_{(k,\ell)} M_k^{\tau-1-\ell}$. Obviously $\varphi_{(k,j)}(t)$ and $a_{(k,\ell)}$ take the same  explicit expressions as presented in Proposition \ref{thm : prop-f_j} ; for the case of a unique eigenvalue. Set $b_{(k,0)}=1$ and $b_{(k,\ell)}=-a_{(k,\ell-1)}$ for $\ell=1,\ldots,\tau$. Hence, we have
 \begin{equation*}
 M_{(k,\tau)} = \sum_{\ell=0}^{\tau} b_{(k,\ell)} M_k^{\tau-\ell} = b_{(k,\tau)}I_{d_k}+\sum_{\ell=0}^{\tau-1} b_{(k,\ell)} M_k^{\tau-\ell}.
 \end{equation*}
 Besides, the $\tau-th$ powers of the matrices $M_k$ may take the form
 \begin{equation*}
 M_k^{\tau} = \left[(M_k-\lambda_k I_{d_k})+\lambda_k I_{d_k}\right]^{\tau}
  =\lambda_k^{\tau}I_{d_k}+\sum_{\eta=1}^{\tau}
  \binom{\eta}{\tau}
  \lambda_k^{\tau-\eta} (M_k-\lambda_k I_{d_k})^{\eta} .
 \end{equation*}
 %%%
 Thus, we obtain
 \begin{equation*}
 M_{(k,\tau)}=\sum_{\ell=0}^{\tau} b_{(k,\ell)} \lambda_k^{\tau-\ell}I_{d_k}+
 \sum_{\ell=0}^{\tau-1} b_{(k,\ell)}\sum_{\eta=1}^{\tau-\ell}
 \binom{\eta}{\tau-\ell}
 \lambda_k^{\tau-\ell-\eta} (M_k-\lambda_k I_{d_k})^{\eta}.
 \end{equation*}
 %%%
 Consequently, we have
 {\small
 \begin{eqnarray*}
 f_{[j]}(I_d-tA) &=& \sum_{k=1}^l \exp(2i \pi j_k/p)
 \left[ \varphi_{(k,0)}(t)+\sum_{\tau=1}^{m_k-1} \varphi_{(k,\tau)}(t)\sum_{\ell=0}^{\tau} b_{(k,\ell)} \lambda_k^{\tau-\ell} \right] I_{(d_k,A)}\\
 &+& \sum_{k=1}^l \exp(2i \pi j_k/p)
 \left[ \sum_{\tau=1}^{m_k-1}\varphi_{(k,\tau)}(t) \Upsilon(k,\tau) (M_k-\lambda_k I_{d_k})_A^{\eta}\right],
 \end{eqnarray*}
 }
 where $\Upsilon(k,\tau) =\Upsilon(A,k,\tau)= \sum_{\ell=0}^{\tau-1} b_{(k,\ell)}\sum_{\eta=1}^{\tau-\ell}\binom{\eta}{\tau-\ell}
 \lambda_k^{\tau-\ell-\eta}$, $I_{(d_k,A)}=\Theta_{d_{1}} \oplus \cdots \oplus
 \Theta_{d_{k-1}} \oplus I_{d_k} \oplus
 \Theta_{d_{k+1}} \oplus \cdots \oplus \Theta_{d_l}$ and $(M_k-\lambda_k I_{d_k})_A=\Theta_{d_{1}} \oplus \cdots \oplus
 \Theta_{d_{k-1}} \oplus (M_k(\lambda_k)-\lambda_k I_{d_k})\oplus
 \Theta_{d_{k+1}} \oplus \cdots \oplus \Theta_{d_l}$.
 It was established in Theorem 2 of \cite{ben-rach-2002} that the $\eta-th$ powers of this matrix are given by
 %%%%%%%%%%%%%%%%
 \begin{equation}\label{eq:(M-k-lambda-k)-A^eta}
   (M_k(\lambda_k)-\lambda_k I_{d_k})_A^{\eta}=
 \prod_{\omega=1,\omega \neq k}^{l}
 \frac{(A-\lambda_{\omega}I_d)^{m_{\omega}}}{(\lambda_k-\lambda_{\omega})
 ^{m_{\omega}}} \sum_{i=0}^{m_k-\eta-1} \alpha_{i,k} (A-\lambda_k I_d)^{\eta+i},
 \end{equation}
 %%%%%%%%%%%%%%%%%%%%%%%%%
 %%%%%%%%%%%%%%%%%%%%%%%%%
 with $\alpha_{0,k}=1$ and
 \begin{equation*}%\label{alpha(i,j)-a(i,j)}
 \alpha_{i,k} = \frac{-1}{a_{0,k}}
 \sum_{\theta=1}^i b_{\theta,k}\alpha_{i-\theta,k}\;
 \mbox{ with }   \; b_{i,k} = \sum_{\Gamma_{i,k}}
 \prod_{t=1,t\neq k}^{l}
 \binom{h_t}{m_t}
 (\lambda_k - \lambda_t)^{m_t-h_t},
 \end{equation*}
 where $b_{i,k}=0$ for $i>m_{1}+\dots+m_{k-1}+m_{k+1}+\dots+m_{s}$ and  $\Gamma_{i,k}$ is the set of $(h_1,\ldots, h_{k-1},h_{k+1},\ldots,h_l) \in \N^{l-1}$ satisfying the relation $h_{1}$+$\dots$+$h_{k-1}+h_{k+1}$+$\dots$+$h_{l}=i$ with $h_t\leq m_t$.  More precisely, in \cite{bent-rach-2008}, an explicit formula of $\alpha_{i,k}$ is provided as follows
 \begin{equation}\label{eq:alpha(i,k)}
   \alpha_{i,k}=(-1)^i \sum_{\Gamma_{i,k}}
 \prod_{t=1,t\neq k}^{l}
 \binom{h_t}{m_t+h_t-1}
 (\lambda_k - \lambda_t)^{-h_t}.
 \end{equation}
 %%%%
 \par
 %%%%
 For the sake of simplicity, in the remainder of this section, it will be denoted by $A$ a nonsingular matrix, in $M_d(\C)$, reduced to its Jordan canonical form. More precisely, whether $A$ admits $l$ eigenvalues $\lambda_1,\ldots,\lambda_l$, $A$ will take the form
 $\bigoplus_{k=1}^l M_k$, where the matrix $M_k=M_k(\lambda_k)$ is the direct sum of all Jordan blocks associated to $\lambda_k$; which size will be denoted by $d_k$ $(\sum_{k=1}^{l}d_k=d)$.
 %%%%
 \par
 The above discussion is summarized in the following proposition,
 %%%%%%%%%%%%
 \begin{Proposition}\label{thm : prop-pth-roots}
   Let $A \in M_d(\C)$ with minimal polynomial $M_A(z)=\prod_{i=1}^{l}(z-\lambda_i)^{m_i}$ $(\sum_{i=1}^{l} m_i \leq d)$ and $t \in \R \backslash \{0\}$ such that $|t|\sigma(A)<1$.
 Then,
 {\small
 \begin{eqnarray*}
 f_{[j]}(I_d-tA) &=& \sum_{k=1}^l \exp(2i \pi j_k/p)
 \left[ \varphi_{(k,0)}(t)+\sum_{\tau=1}^{m_k-1} \varphi_{(k,\tau)}(t)\sum_{\ell=0}^{\tau} b_{(k,\ell)} \lambda_k^{\tau-\ell} \right] I_{(d_k,A)}\\
 &+& \sum_{k=1}^l \exp(2i \pi j_k/p)
 \left[ \sum_{\tau=1}^{m_k-1} \varphi_{(k,\tau)}(t) \Upsilon(k,\tau)  (M_k-\lambda_k I_{d_k})_A^{\eta}\right],
 \end{eqnarray*}
 }
 where $j \equiv (j_1,\ldots,j_l)$ $(j_k \in R(p))$ and the $\eta-th$ powers of $(M_k-\lambda_k I_{d_k})_A$ are given by (\ref{eq:(M-k-lambda-k)-A^eta})-(\ref{eq:alpha(i,k)}).
 \end{Proposition}
 %%%%
 Explicit formula for all primary matrix $pth$ roots may be obtained using the Lagrange-Sylverster interpolation polynomial \cite{gantmacher-1960}. However, as far as we know our formula presented in Proposition \ref{thm : prop-pth-roots} is not known under this form in the literature.
 %%%%
 \par
 The following two corollaries are immediately derived from the Proposition \ref{thm : prop-pth-roots}.
 %%%
 \begin{Corollary}\label{thm : coro-2eig}
   Consider $A$ with minimal polynomial $M_A(z)=(z-\lambda_1)^{m_1}(z-\lambda_2)^{m_2}$ and $t$ a non-zero real number, such that $|t|<\frac{1}{\max \{|\lambda_1|,|\lambda_2|\}}$. Then,
 {\small
 \begin{eqnarray*}
 f_{[j]}(I_d-tA) &=& \sum_{k=1}^2 \exp(2i \pi j_k/p)
 \left[ \varphi_{(k,0)}(t)+\sum_{\tau=1}^{m_k-1} \varphi_{(k,\tau)}(t)\sum_{\ell=0}^{\tau} b_{(k,\ell)} \lambda_k^{\tau-\ell} \right] I_{(d_k,A)}\\
 &+& \sum_{k=1}^2 \exp(2i \pi j_k/p)
 \left[ \sum_{\tau=1}^{m_k-1} \varphi_{(k,\tau)}(t) \Upsilon(k,\tau) (M_k-\lambda_k I_{d_k})_A^{\eta}\right],
 \end{eqnarray*}
 }
 where $[j] \equiv (j_1,j_2)$ $(j_k \in R(p))$ and the $\varphi_{(k,\tau)}$ are as in Proposition \ref{thm : prop-f_j}.
 \end{Corollary}
 %%%
 \begin{Example}
  Compute the square roots of the $4 \times 4$ matrix given by,
  \begin{equation*}
  B=\left(\begin{array}{cccc}
  2/3 & -1 & 0 & 0 \\
  0 & 2/3 & 0 & 0 \\
  0 & 0 & 1/3 & -1 \\
  0 & 0 & 0 & 1/3
  \end{array}
  \right).
  \end{equation*}
  %%%
  For $A=I_4-B$, we have
  $M_A(z)=(z-\lambda_1)^2(z-\lambda_2)^2$ with $\lambda_1=1/3$ and $\lambda_2=2/3$. Therefore, according to Corollary \ref{thm : coro-2eig}, the matrix square roots of $B=I_4-A$ are given as follows,
  \begin{eqnarray*}
   f_{[j]}(B) &=& \sum_{k=1}^2 \exp(i \pi j_k)
 \left[ \varphi_{(k,0)}(1)+\varphi_{(k,1)}(1)\sum_{\ell=0}^{1}
  b_{(k,\ell)} \lambda_k^{1-\ell} \right] I_{(2,A)} \\
  &+& \sum_{k=1}^2 \exp(i \pi j_k)
  \varphi_{(k,1)}(1) b_{(k,0)} (M_k-\lambda_k I_{2})_A %\\
  \end{eqnarray*}
 where $[j] \equiv (j_1,j_2) \in (\{0,1\})^2$, $b_{(1,0)}=b_{(2,0)}=1$, $b_{(1,1)}=-2/3$, $b_{(2,1)}=-4/3$,
 $\varphi_{(1,0)}(1)=\sqrt{\frac{2}{3}}-\frac{1}{6}\sqrt{\frac{3}{2}}$,
 $\varphi_{(1,1)}(1)=-\frac{1}{2}\sqrt{\frac{3}{2}}$, $\varphi_{(2,0)}(1)=0$ and $\varphi_{(2,1)}(1)=-\frac{\sqrt{3}}{2}$.
 Consequently the $4\times 4$ square roots of the matrix $B$ correspond to the following matrices,
  \begin{eqnarray*}
  && \left(\begin{array}{cccc}
  \sqrt{\frac{2}{3}} & -\frac{1}{2}\sqrt{\frac{3}{2}} & 0 & 0 \\
  0 & \sqrt{\frac{2}{3}} & 0 & 0 \\
  0 & 0 & \frac{\sqrt{3}}{3} & -\frac{\sqrt{3}}{2} \\
  0 & 0 & 0 & \frac{\sqrt{3}}{3}
  \end{array}
  \right); \;\
  \left(\begin{array}{cccc}
  -\sqrt{\frac{2}{3}} & \frac{1}{2}\sqrt{\frac{3}{2}} & 0 & 0 \\
  0 & -\sqrt{\frac{2}{3}} & 0 & 0 \\
  0 & 0 & \frac{\sqrt{3}}{3} & -\frac{\sqrt{3}}{2} \\
  0 & 0 & 0 & \frac{\sqrt{3}}{3}
  \end{array}
  \right);\\
  && \left(\begin{array}{cccc}
  \sqrt{\frac{2}{3}} & -\frac{1}{2}\sqrt{\frac{3}{2}} & 0 & 0 \\
  0 & \sqrt{\frac{2}{3}} & 0 & 0 \\
  0 & 0 & -\frac{\sqrt{3}}{3} & \frac{\sqrt{3}}{2} \\
  0 & 0 & 0 & -\frac{\sqrt{3}}{3}
  \end{array}
  \right)  ; \;\
  \left(\begin{array}{cccc}
  -\sqrt{\frac{2}{3}} & \frac{1}{2}\sqrt{\frac{3}{2}} & 0 & 0 \\
  0 & -\sqrt{\frac{2}{3}} & 0 & 0 \\
  0 & 0 & -\frac{\sqrt{3}}{3} & \frac{\sqrt{3}}{2} \\
  0 & 0 & 0 & -\frac{\sqrt{3}}{3}
  \end{array}
  \right).
  \end{eqnarray*}
  %%%%
  \end{Example}
  %%%%%%%%
  %%%%%%%
 \begin{Corollary}\label{thm : coro-3eig}
 Consider $A$ with minimal polynomial $M_A(z)=(z-\lambda_1)^{m_1}(z-\lambda_2)^{m_2}(z-\lambda_3)^{m_3}$ and $t$ a non-zero real number, such that $|t|<\frac{1}{\max \{|\lambda_1|,|\lambda_2|,|\lambda_3|\}}$. Then,
 {\small
 \begin{eqnarray*}
 f_{[j]}(I_d-tA) &=& \sum_{k=1}^3 \exp(2i \pi j_k/p)
 \left[ \varphi_{(k,0)}(t)+\sum_{\tau=1}^{m_k-1} \varphi_{(k,\tau)}(t)\sum_{\ell=0}^{\tau} b_{(k,\ell)} \lambda_k^{\tau-\ell} \right] I_{(d_k,A)}\\
 &+& \sum_{k=1}^3 \exp(2i \pi j_k/p)
 \left[ \sum_{\tau=1}^{m_k-1} \varphi_{(k,\tau)}(t) \Upsilon(k,\tau) (M_k-\lambda_k I_{d_k})_A^{\eta}\right],
 \end{eqnarray*}
 }
 where $[j] \equiv (j_1,j_2,j_3)$ $(j_k \in R(p))$ and the $\varphi_{(k,\tau)}$ are as in Proposition \ref{thm : prop-f_j}.
 \end{Corollary}
 %%%%
  \begin{Example}
  Compute the square roots of the $5 \times 5$ matrix,
  \begin{equation*}
  B=\left(\begin{array}{ccccc}
  1/2 & -1 & 0 & 0 & 0\\
  0 & 1/2 & 0 & 0 & 0\\
  0 & 0 & 2/3 & -1 & 0 \\
  0 & 0 & 0 & 2/3 & 0 \\
  0 & 0 & 0 & 0 & 3/4
  \end{array}
  \right).
  \end{equation*}
  %%%
  For $A=I_5-B$, we have
  $M_A(z)=(z-\lambda_1)^2(z-\lambda_2)^2 (z-\lambda_3)$ with $\lambda_1=1/2$, $\lambda_2=1/3$ and $\lambda_3=1/4$. Therefore, according to Corollary \ref{thm : coro-3eig} the matrix square roots of $B=I_5-A$ are given by,
  \begin{eqnarray*}
  f_{[j]}(B) &=& \sum_{k=1}^2 \exp(i \pi j_k)
 \left[ \varphi_{(k,0)}(1)+\varphi_{(k,1)}(1)\sum_{\ell=0}^{1}
  b_{(k,\ell)} \lambda_k^{1-\ell} \right] I_{(d_k,A)} \\
  &+& \sum_{k=1}^2 \exp(i \pi j_k)
  \varphi_{(k,1)}(1) b_{(k,0)} (M_k-\lambda_k I_{2})_A+
  \exp(i \pi j_3) \varphi_{(3,0)}(1)I_{(d_3,A)},
  \end{eqnarray*}
 where $[j] \equiv (j_1,j_2,j_3) \in (\{0,1\})^3$, $b_{(1,0)}=b_{(2,0)}=1$, $b_{(1,1)}=-1$, $b_{(2,1)}=-2/3$,
 $\varphi_{(1,0)}(1)=\frac{\sqrt{2}}{4}$,
 $\varphi_{(1,1)}(1)=-\frac{\sqrt{2}}{2}$, $\varphi_{(2,0)}(1)=\frac{3}{4}\sqrt{\frac{2}{3}}$, $\varphi_{(2,1)}(1)=-\frac{1}{2}\sqrt{\frac{3}{2}}$ and $\varphi_{(3,0)}(1)=\sqrt{\frac{3}{4}}$.
 Consequently, the $8$ square roots of the matrix $B$, computed using all primary matrix functions of $B$, are presented by the matrices,
 {\small
 \begin{eqnarray*}
 f_{[j]}(B)=\left(\begin{array}{ccccc}
  \exp(i \pi j_1)\frac{\sqrt{2}}{2} & -\exp(i \pi j_1)\frac{\sqrt{2}}{2} & 0 & 0 & 0\\
  0 & \exp(i \pi j_1)\frac{\sqrt{2}}{2} & 0 & 0 & 0\\
  0 & 0 & \exp(i \pi j_2)\sqrt{\frac{2}{3}} &
  -\frac{\exp(i \pi j_2)}{2}\sqrt{\frac{3}{2}} & 0 \\
  0 & 0 & 0 & \exp(i \pi j_2)\sqrt{\frac{2}{3}} & 0 \\
  0 & 0 & 0 & 0 & \exp(i \pi j_3)\sqrt{\frac{3}{4}}
  \end{array}
  \right),
 \end{eqnarray*}
   }
  where $[j] \equiv (j_1,j_2,j_3) \in (\{0,1\})^3$.
   \end{Example}

%{\bf Received: Month xx, 20xx}

\end{document}